\documentclass{article}
\usepackage[utf8]{inputenc}

\title{Goldbach's Racetrack}
\author{René van der Vegt }
\date{August 2020}

\usepackage{graphicx}
\usepackage{amsmath}
\usepackage{relsize}%
\usepackage{float}
\usepackage{natbib}
\begin{document}

\maketitle

\section*{Abstract}

The Goldbach conjecture that "every even integer is the sum of two primes" has been open since 1742. This paper details a road map to a proof of Goldbach’s conjecture based on a function that estimates the number of Goldbach pairs. It is shown that the estimated of the number of Goldbach pairs increases as the even integer increases and that the error term is small. Further detailing the functions, assumptions and error term may finally lead to a proof of Goldbach's conjecture.

\section{Introduction}

Goldbach’s conjecture resulted from a letter exchange between Euler and Goldbach. Euler stated on the 30th of June 1742 that “every even integer is the sum of two primes” \citep{Goldbach}. The conjecture has been numerically verified up to $4 \times 10^{18}$ \citep{Silva:2013}.

Let $2n$ be an even number. A Goldbach Pair (GP) of 2n consists of two prime numbers whose sum is $2n$. The approach to proof the Goldbach conjecture is to define a function that estimates the number of Goldbach pairs for a given $2n$. If you can show that value of this function is forever increasing and that the remaining error term between the estimate and the actual number of GP is small and bounded, than you have proven the Goldbach conjecture. 

\begin{figure}[H]
    \centering
    \includegraphics[width=\linewidth]{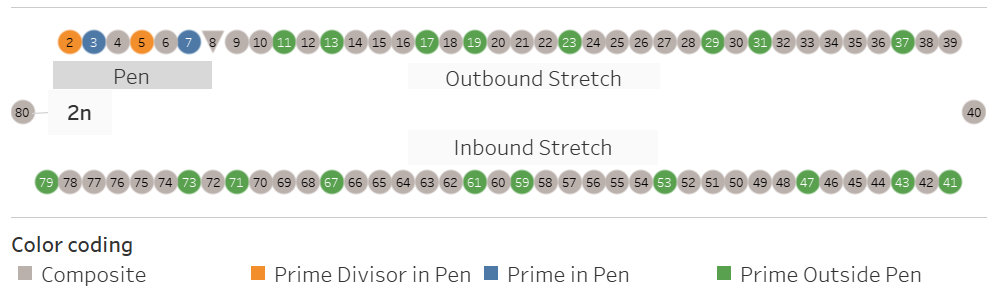}
    \caption{Goldbach's Racetrack for 2n=80}
    \label{fig:racetrack}
\end{figure}

As a starting point, we construct a Goldbach racetrack for estimating the number of Goldbach Pairs as schematically laid out in Fig. \ref{fig:racetrack}. 
For an even number $2n$ there are initially $n$ pairs that have a sum $2n$. For a pair to be a GP, both members of the pair have to be prime. To determine the GP a method similar to the sieve of Eratosthenes is used to filter the pairs. To eliminate all the composite numbers (grey) below $2n$ you only have to use the sieve for primes up to $\sqrt{2n}$, the Pen. All the primes in the Pen go around the racetrack from low to high. All multiples of the prime are composite numbers. The associated pair has at least one composite number and is eliminated. Once all the primes in the Pen have gone round the racetrack, you have only primes remaining up to $2n$ and the remaining pairs are the GP for $2n$. 

 As 2 goes around the racetrack half of the pairs are eliminated, i.e. all the pairs that consist of even numbers, and half remain standing. As 3 goes around then 1/3 of the remaining pairs are eliminated on the outbound stretch and an additional 1/3 on the inbound stretch (assuming that 3 is not a divisor of $2n$). 
 For a prime $p>2$ in the Pen that goes round the racetrack and is not a divisor of $2n$ (blue), the fraction of pairs that is eliminated is $\frac{2}{p}$; $\frac{1}{p}$ on the outbound stretch and $\frac{1}{p}$ on the inbound stretch. The fraction that remains standing is $1-\frac{2}{p} = \frac{(p-2)}{p}$. If the prime in the Pen that goes round the racetrack is a divisor of $2n$ (orange), than on the inbound part of the racetrack, the prime eliminates the same pairs as on the outbound stretch, i.e. the fraction that remains standing is $\frac{p-1}{p}$. In the case for $2n=80$, there are 4 primes in the Pen, 2 and 5 are divisors of 80 and 3 and 7 are not. The fraction of pairs remaining standing are $\frac{(2-1)}{2}\ \frac{(3-2)}{3}\ \frac{(5-1)}{5}\ \frac{(7-2)}{7} = \frac{2}{21}$ The estimated number of GP for 2n=80 would be $\frac{2}{21}\  40 \approx 4$. 

Generalizing gives that the estimated number of GP for $2n$ (EGP(2n)) is 

\begin{equation}
\begin{split}
\ & EGP(2n) =\frac{n}{2} \ \prod_{p=3}^{p<\sqrt{2n}} \frac{f(2n,p)}{p} \\ 
 & \text{   where   }  f(2n,p) \begin{cases}
      p-1,&\text{if}\ mod(2n,p)=0 \\
      p-2,&\text{if}\ mod(2n,p) \ne 0\\
    \end{cases}
    \label{eq:EGP}
\end{split}
\end{equation}

\begin{figure}[H]
    \centering
    \includegraphics[width=\linewidth]{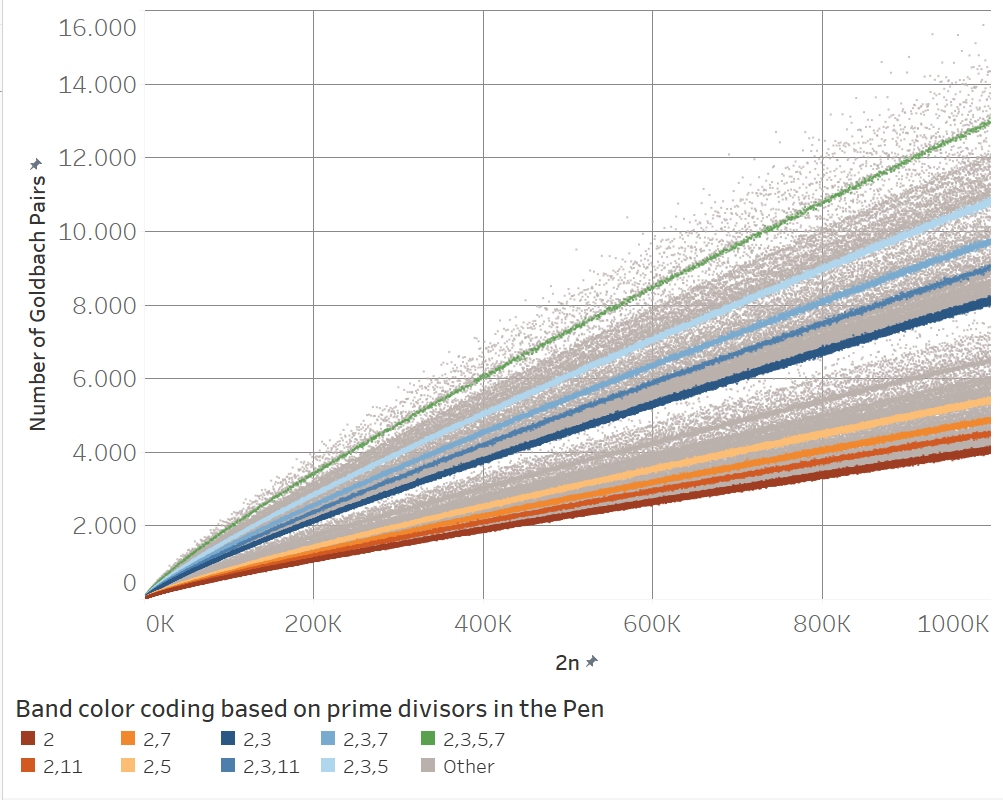}
    \caption{Goldbach Pair in Band's}
    \label{fig:bands}
\end{figure}

The estimate of the GP depends on whether or not the primes in the Pen are divisors of $2n$. If the estimate is accurate you would expect that the actual number of GP for $2n$ also varies in the same manner. Fig. $\ref{fig:bands}$ shows the actual number of GP for all $2n$ up to 1 Million. The numbers $2n$ are grouped together in bands (B) based on the prime divisors in the Pen, e.g. $B_2$ denotes the group of $2n$ with only 2 as prime divisor in the Pen. The colors in the Fig. $\ref{fig:bands}$ denote 9 specific bands, all 2n that belong to other bands are grey. 

The distribution of the actual number of GP follows the pattern that you would expect from Eq. (\ref{eq:EGP}). The Band $B_2$ has the lowest number of GP. As more primes in the Pen are divisors of $2n$, the higher the actual number of GP. If you compare $B_{2,3}$ (dark blue) and $B_{2}$ (dark red) then both have the same product terms in Eq. (\ref{eq:EGP}) except for $p=3$. For $p=3$, $B_{2,3}$ has the term $2/3$ while $B_{2}$ has the term $1/3$. Therefore the estimate for $B_{2,3}$ (dark blue) is twice that for $B_2$ (dark red). Fig. $\ref{fig:bands}$ shows that the ratio of the actual number of GP for $B_2$ and $B_{2,3}$ is the same ratio as from Eq. (\ref{eq:EGP}).

\begin{figure}[H]
    \centering
    \includegraphics[width=\linewidth]{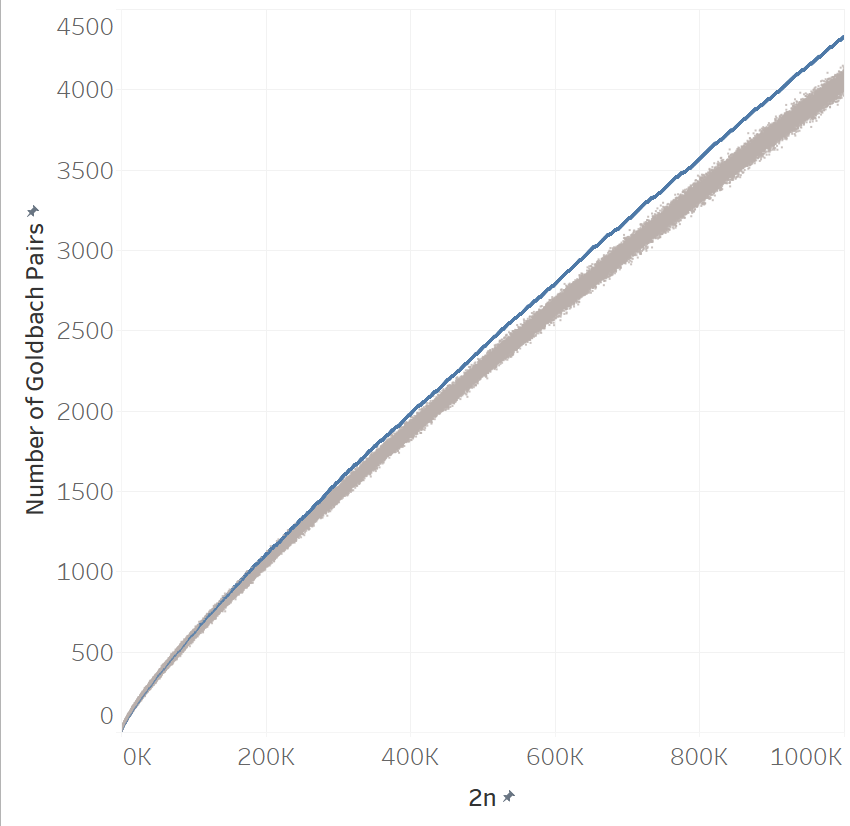}
    \caption{Estimated Goldbach Pairs (blue) vs Actual Goldbach Pairs (grey) for $B_2$}
    \label{fig:EGP and GP for B2}
\end{figure}

 Fig. \ref{fig:EGP and GP for B2} compares the EGP (blue) and actual GP (grey) for band $B_2$. EGP slightly overestimates the number of GP. 

It seems that Eq. (\ref{eq:EGP}) gives a good estimate of the number of GP for a given $2n$. It explains the pattern found in the actual number of GP but slightly overestimates the number of GP.

\newpage
\section{Elimination probability function}

The approximation used in the initial estimate of the number of GP, Eq. (\ref{eq:EGP}), is that the average probability of eliminating a number and the associated pair is equal to $1/p$. 
If $p$ goes around the racetrack then up to $p^2$ the probability of eliminating a pair is 0, as all the multiples of $p$ have already been eliminated by smaller primes. Therefore if $p$ is just smaller than $\sqrt(2n)$, the average probability of eliminating a pair as $p$ goes around the racetrack is much smaller than $1/p$.

To improve the estimate of the number of GP we need to define a function for the probability of a multiple of $p$ as a function of a number $x$. On the racetrack the primes in the Pen go round from low to high. Thus by the time $p$ goes round the racetrack all the multiples of $p$ with divisors smaller then $p$ have already been eliminated. When $p$ goes round the racetrack it can only eliminate a number $x$ if $x$ has only divisors larger than or equal to $p$. Therefore the probability of an elimination by $p$ for a number $x$ smaller than $p^2$ is 0. For numbers larger than $p^2$, an elimination only occurs if that number does not have divisors smaller then $p$ and is divisible by $p$. If it does have smaller divisors, the number and associated pair, will already have been eliminated by the smaller prime.  

You can split all numbers in two groups based on a prime $p$. All numbers whose prime factors are higher or equal to $p$ ($H_p$) and the rest ($L_p$). As $p$ goes around the racetrack it can only eliminate numbers that are in the $H_p$ group. The fraction $F_H$ of all numbers for a given $p$ is 

\begin{equation}
\ F_H(p)= \prod_{k=2}^{k<p} \frac{k-1}{k}   \ k\  is\  prime  
  \label{eq:FH}
\end{equation}

\begin{figure}[H]
    \centering
    \includegraphics[width=\linewidth]{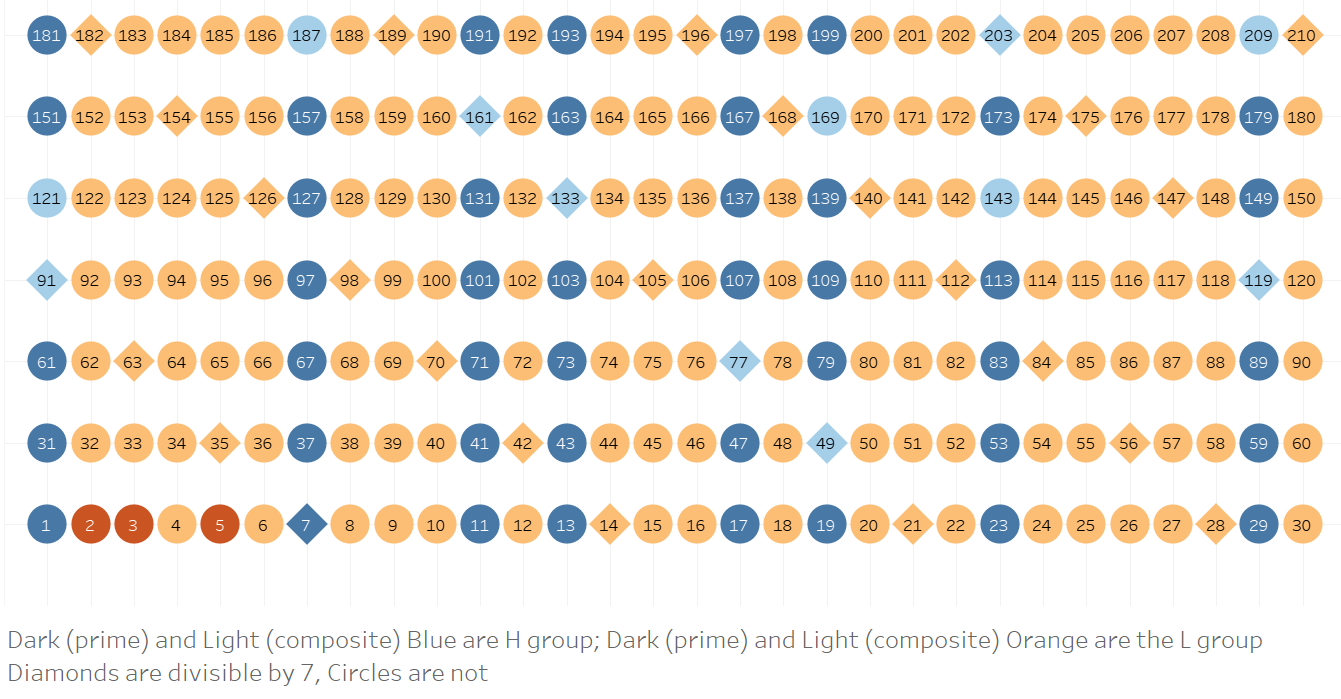}
    \caption{H and L group for p=7}
    \label{fig:H and L}
\end{figure}

Fig. \ref{fig:H and L} illustrates the distribution of numbers in $H_7$ and $L_7$ groups for $p=7$. The numbers are arranged in 7 rows and each row contains 30 numbers. 30 is the primorial of 5 (5\#), i.e. the primorial of the biggest prime smaller than 7. Numbers are colour coded as dark blue (prime, $\in H$), light blue (composite, $\in H$) and dark orange (prime, $\in L$) and light orange (composite, $\in L$). The fraction of $H_7$ on the first row according to Eq. (\ref{eq:FH}) is 1/2 $\times$ 2/3 $\times$ 4/5 = 8/30, i.e there are 8 blue numbers on the first row. Note that 1 is considered $\in H_p$.

As 7 and 30 are co-prime and all the primes smaller than 7 are divisors of 30, the numbers in a column belong to the same group and have the same color either blue or orange. The fraction of H ($F_H(7)$) in each row is the same and equal to Eq. (\ref{eq:FH}). As 7 and 30 are co-prime and there are 7 rows, there is exactly one number which is divisible by 7 in each column. A number divisible by 7 is shown as a diamond in Fig. \ref{fig:H and L}. 

 The same principal as illustrated in Fig. \ref{fig:H and L} applies for any prime $p$.
 Each column belongs to either $H_p$ or $L_p$ group and in each column there is exactly one number divisible by $p$. Therefore the average probability that $x$ is divisible by $p$ for all $\{{x\in H_p\ |\  x\leq p\#\}}$ is $1/p$ 

\begin{equation}
\ \overline{\Pr\big(mod(x,p)=0\big)} = \frac{1}{p}  \ \ \ \   {\{x\in H_p\ |\  x\leq p\#\}} 
  \label{eq:average Pr H}
\end{equation}

For numbers larger than $p\#$ the principles as illustrated in Fig. \ref{fig:H and L} can be extended. As more rows are added, they stay either blue or orange. If a number $x \leq 210$ is divisible by 7 than $x + 210$ will also be divisible by 7. 
The probability for number $x$ to be divisible by $p$ for a number $x > p\#$ is equal to the probability of $mod(x,p\#)$ being divisible by $p$. In other words the probability distribution is cyclical with a period of $p\#$ 

\begin{equation}
\ \Pr(mod(x,p)=0) \ =\  \Pr(mod(mod(x,p\#),p)=0  \ \ \{{x\in H_p\ |\  x>p\#\}} 
  \label{eq:cyclical probability}
\end{equation}

\newpage
\section{Elimination probability per factor multiple}

The previous section established some general properties of the probability function that $x \in H_p$ is divisible by $p$. Up to $p^2$ the probability is 0 and we need to quantify the probability function between $p^2$ and $p\#$. If $x$ is divisible by $p$, than $x$ can be a two, three or higher factor multiple of $p$. As $x \in H_p$, all the factors of $x$ have to be larger or equal than $p$.

To quantify the overall probability that $x \in H_p$ is divisible by $p$, we need to quantify the probability for each factor multiple and sum over all potential factors to find the overall probability.

\subsection{Two factor multiples}

Let $x$ be a two-factor multiple $p \times r$, where $p\leq r$ and $x\leq p\#$. According to the Prime Number Theorem \citep{Hadamard1896},\citep{Poussin1896}, the probability that $r$ is prime is $\mathlarger{\frac{1}{ln(r)} = \frac{1}{\ln(\frac{x}{p})}}$. 

The probability that $x$ is a two-factor multiple of p equals the probability of $x$ being divisible by $p$ $(1/p)$ times the probability that $r$ is prime :  $\mathlarger{\frac{1}{p} \frac{1}{\ln(\frac{x}{p})}}$. This probability is relative to all numbers but we are interested in the probability for $x \in H_p$. The absolute probability of a two-factor multiple of $p$ for $x\in H_p$ is 

\begin{equation}
\ Pr_{2f}(x,p)= \mathlarger { \frac{1}{p} \frac{1}{\ln(\frac{x}{p})} \frac{1}{F_{H_p}} }\ \ \   \{x\in{H_p\ |\  p^2\leq x\leq p\#\}}
     \label{eq:Pr 2factor multiples}
\end{equation}

Instead of quantifying the absolute probability as in Eq. (\ref{eq:Pr 2factor multiples}), it is also possible to quantify the probability in relative terms to the average probability of $1/p$ assumed in Eq. (\ref{eq:EGP}). Eq. (\ref{eq:average Pr H}) establishes that the average probability over the range 1...$p\#$ is $1/p$. 

\begin{equation}
\ Pr_{2f}(x,p)= RPF_{2f}(x,p) \times  \frac{1}{p}    \ \ \ \ {\{x\in H_p\ |\  p^2\leq x\leq p\#\}}
     \label{eq:RPF}
\end{equation}

To convert the absolute probability to the Relative Probability Factor (RPF) we need to multiple by $p$.  

As we want a proof of the Goldbach conjecture for all numbers, we want to establish the probability as $p \to \infty$. We can do this by applying Merten's third Theorem \citep{Mertens1874} , which states that

\begin{equation}
\ \lim_{n\to\infty} \frac{1}{\ln(n)} \prod_{p\leq n} (1-\frac{1}{p})  = e^{-\gamma}
   \label{eq:mertens 3rd}
\end{equation}

\noindent
and substituting Eq. (\ref{eq:FH}) and Eq. (\ref{eq:Pr 2factor multiples}) in Eq. (\ref{eq:RPF}). This gives the $RPF_{2f}$ as $p \to \infty$

\begin{equation}
\ RPF_{2f}(x,p)=      \mathlarger { \frac{e^\gamma}{\log_p(\frac{x}{p})} } \ \   \{{x\in H_p\ |\  p^2\leq x\leq p\#\}}
     \label{eq:RPF for 2 factor multiples }
\end{equation}

\begin{figure}[H]
    \centering
    \includegraphics[width=\linewidth]{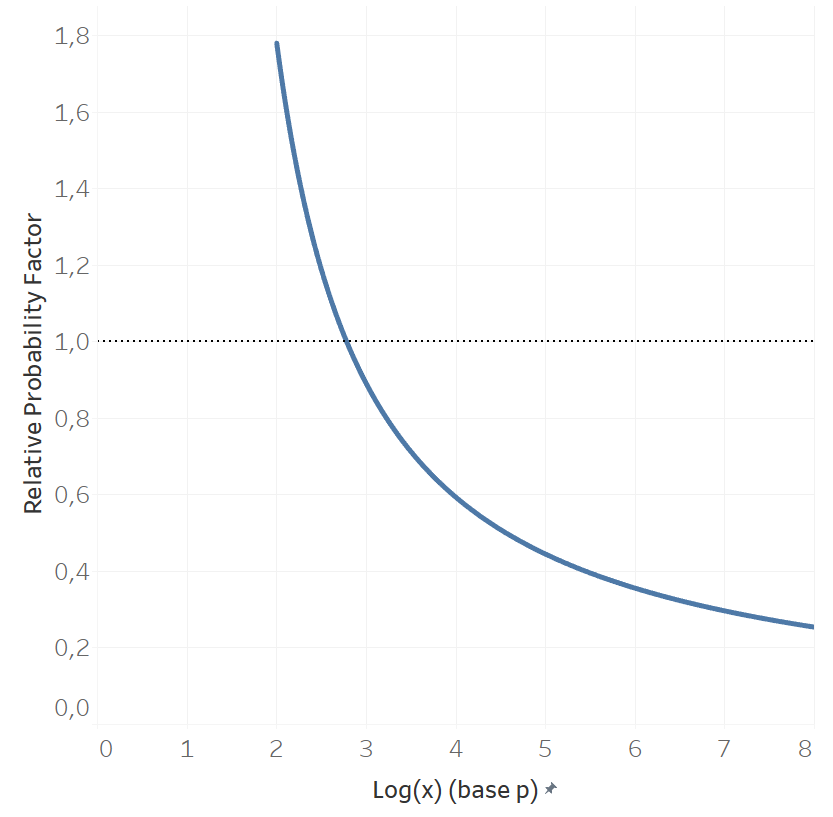}
    \caption{RPF for 2 factor multiples of p, as $p \to \infty,\  x \in H_p$}
    \label{fig:RPF 2factor multiple}
\end{figure}

Fig. \ref{fig:RPF 2factor multiple} shows the value of $RPF_{2f}$.

The $RPF_{2f}$ of finding a number $x \in H_p$ that is a two factor multiple of p is 0 up to $p^2$, then shoots up to $e^\gamma$ and then declines with $1/{log_p (x/p)}$.

\subsection{Higher factor multiples}

Three factor multiples with all factors larger or equal than p start from $p^3$. All three factor multiples can be written as $p \times r \times s$. where $p\leq r \leq s$. The probability that $x$ is a three-factor multiple of $p$, given that $p$ and $r$ are prime equals the probability of $x$ being divisible by $p$, $(1/p)$ and divisible by $r$ $(1/r)$ times the probability that s is prime:  $\mathlarger{\frac{1}{pr} \frac{1}{\ln(\frac{x}{pr})}}$. The absolute probability for a three factor multiple is.

\begin{equation}
\ Pr_{3f}(x,p,r)=  \mathlarger { \frac{1}{p r} \frac{1}{\ln(\frac{x}{p r})} \frac{1}{F_{H_p}} } \ \ \   \{{x\in H_p\ |\  p^3\leq x\leq p\#\}}
    \label{eq:Pr 3factor multiples}
\end{equation}

Similarly applying Merten's Third Theorem for a three-factor multiple $(p,r,s)$, the limit of the relative probability factor as $p\to\infty$ at $pr^2$ is $\mathlarger{ \frac{e^\gamma}{r} } $.
The $RPF_{3f}$ for a number $x\in H_p$ and a three-factor multiple $(p,r,s)$ is

\begin{equation}
\ RPF_{3f}(x,p,r)=  \mathlarger{ \mathlarger  { \frac{\frac{e^\gamma}{r}}{\log_p(\frac{x}{p r})} } } \ \ \ \{{x \in H_p\ |\  p^3\leq x\leq p\#\}}
    \label{eq:RPF for 3 factor multiples }
\end{equation}

The shape is similar to the shape of the $RPF_{2f}$ in Fig. \ref{fig:RPF 2factor multiple}. The $RPF_{3f}$ is 0 up to $pr^2$ and then jumps to $\mathlarger{ \frac{e^\gamma}{r} }$. 

There are many possible value of $r$, each with its own $RPF_{3f}$ curve. Each curve starts at $pr^2$ and has a similar shape as in Fig. \ref{fig:RPF 2factor multiple} . As $p\leq r \leq s$ than $r\leq \sqrt \frac{x}{p}$. The sum for all possible three-factor multiple of $p$ is given by summing the $RPF_{3f}$  for each $r$.

\begin{equation}
\ TRPF_{3f}(x,p)= \mathlarger { \mathlarger {\sum_{r=p}^{r=\sqrt \frac{x}{p}} \frac{\frac{e^\gamma}{r}}{\log_p(\frac{x}{pr})} } } \ \ \   \{{x\in H_p\ |\  p^3\leq x\leq p\#\}}
    \label{eq:TRPF for 3 factor multiples}
\end{equation}

The same approach can be applied to higher factor multiples. Four factor multiples start at $p^4$. For four factor multiples $(p,r,s,t)$ where $p\le r \le s \le t$, the total relative probability factor for $x\in H$ is 

\begin{equation}
\ TRPF_{4f}(x,p)= \mathlarger { \sum_{r\ =\ p}^{r\ = \mathlarger{\sqrt[3]{\frac{x}{p}}}}  \ \ \sum_{s\ =\ r}^{s\  =\mathlarger{\sqrt \frac{x}{pr}} }}  \mathlarger{ \mathlarger {\frac{\frac{e^\gamma}{rs}}{\log_p(\frac{x}{prs})} } }  \ \ \  \{{x\in H_p\ |\  p^4\leq x\leq p\#\}}
    \label{eq:TRPF for 4 factor multiples}
\end{equation}

\begin{figure}[H]
    \centering
    \includegraphics[width=\linewidth]{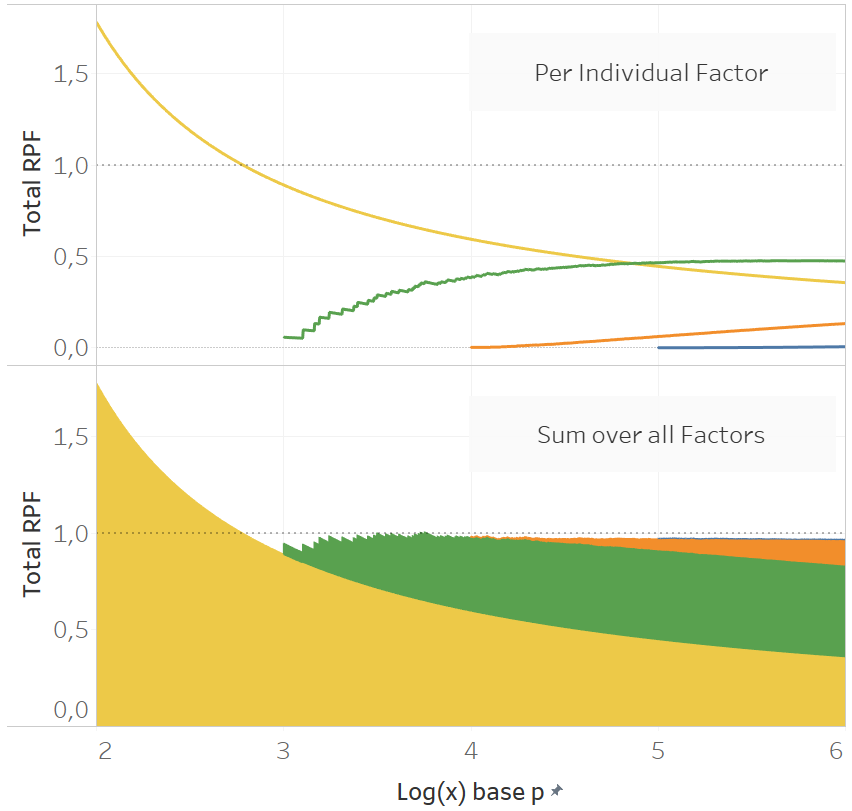}
    \caption{TRPF for 2-5 factor multiples for p=31}
    \label{fig:RPF actual 31}
\end{figure}

Fig. \ref{fig:RPF actual 31} shows the $TRPF$ for 2-5 factor multiples for $p=31$. Note that for two factor multiples the $TRPF$ is the same as $RPF$ as there is only 1 curve.
The top part shows the $TRPF$ for each factor, the bottom part the sum over all factors. The two factor $TRPF_{2f}$ Eq. (\ref{eq:RPF for 2 factor multiples }, yellow) starts at $p^2$ and declines. The three factor $TRPF_{3f}$ Eq. (\ref{eq:TRPF for 3 factor multiples}, green) starts at $p^3$ and for each $r$ there is a curve. Similar for four factor multiples Eq. (\ref{eq:TRPF for 4 factor multiples}, orange) and five factor multiples (blue). The sum of all factors seems to level out just below 1, which is what you expect as the average TRPF over all factors up to $p\#$ is one, according to Eq. (\ref{eq:average Pr H}).

\begin{figure}[H]
    \centering
    \includegraphics[width=\linewidth]{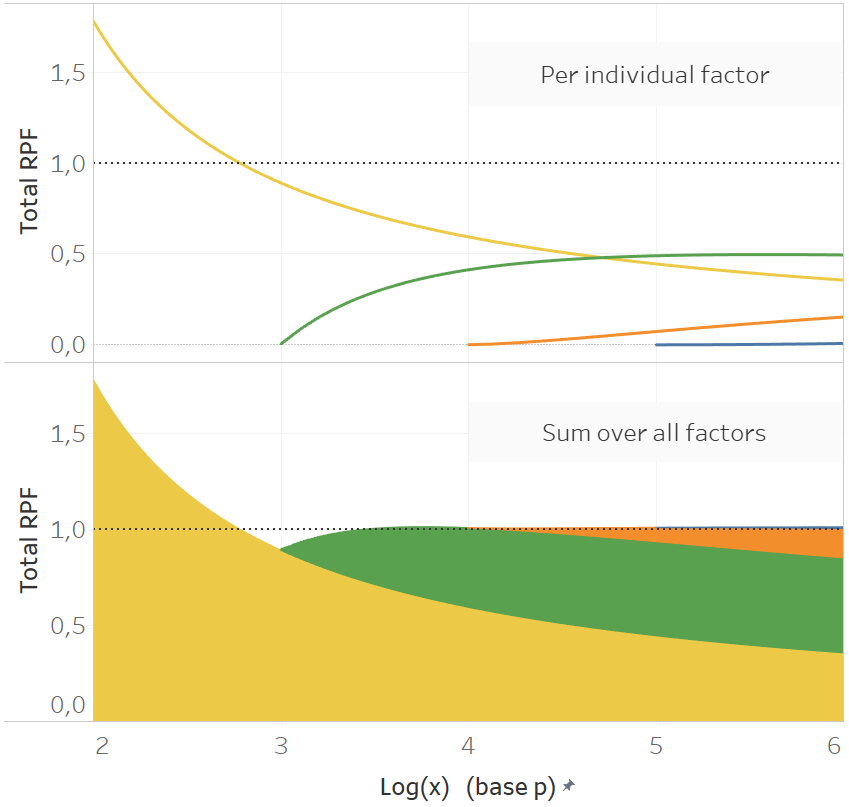}
    \caption{Total RPF for 2-5 factor multiples based on simulated prime distribution}
    \label{fig:RPF ln distribution}
\end{figure}

To emulate what the TRPF curves are for very large primes, the equation for TRPF (Eq. \ref{eq:RPF for 2 factor multiples }, \ref{eq:TRPF for 3 factor multiples} and \ref{eq:TRPF for 4 factor multiples}) are used. In stead of using the actual primes larger than $p$, a simulated prime series is used where the next prime equals $p+\ln(p)$. The approach ensures that the actual numbers used in the emulation remain small enough for computing. The resulting TRPF are shown in Fig. \ref{fig:RPF ln distribution}.
In general it shows the same distribution but it is much smoother due to the fact that many more curves are calculated. Also it levels out just around 1. You would expect the sum of all factor TRPF to level out at 1, as the average TRPF up to $p\#$ is 1, according to Eq. (\ref{eq:average Pr H}).  

\newpage
\section{Improving the GP estimate using TRPF}

The TRPF can be used to improve the estimate of the number of GP for a given $2n$. Fig. \ref{fig:RPF ln distribution} shows that for $\log_p(x) \geq 4$ the TRPF is approximately 1. Up to $\log_p(x) = 2$ the TRPF is 0, between $\log_p(x) = 2$ and $\log_p(x) = 4$ we can estimate the average TRPF using Eq. (\ref{eq:RPF for 2 factor multiples }) and Eq (\ref{eq:TRPF for 3 factor multiples}). 

For each $p$ we want to know the average TRPF up to 2n. With this we can improve the initial estimate that the chance of eliminating a pair is always $1/p$. The average TRPF up to 2n is the integral of the TRPF function divided by 2n.

\noindent
For 2 factor multiples the integral of the TRPF is.

\begin{equation}
\begin{split}
\ I_{2f}(2n,p) & =\mathlarger{ \int_{p^2}^{2n}  \frac{e^\gamma}{\log_p(\frac{x}{p})} } dx \\
            & = e^\gamma\ \ln(p) \  p \  Li \left(\frac{x}{p} \right) \mathlarger{\big|_{p^2}^{2n}} \\
\label{eq: integral 2f}
\end{split}
\end{equation}

\noindent
For 3 factor multiples the integral of the TRPF is.

\begin{equation}
\begin{split}
\ I_{3f}(2n,p) &= e^\gamma\ \ln(p) \left[\left(Li\big(\frac{x}{p}\big)-Li\big(\frac{x}{p^2}\big) \right) \right.  \\
  &\left. - x\bigg( ln(ln(\frac{x}{p}))+ln(2)-ln(ln(\frac{x}{p^2}))\bigg)\right] \mathlarger{\bigg|_{p^3}^{2n}} \\
\label{eq: integral 3f}
\end{split}
\end{equation}

\noindent
We can now define the function $\alpha (2n,p)$ as the average TRPF over all multiples for a given integer $2n$ and a prime $p$ ($p \leq \sqrt{2n}$) (\ref{eq: alpha})

\begin{equation}
\ \alpha(2n,p)= \begin{cases}
      1 &\text{for}\ 4 \leq log_p(2n) \\
      \mathlarger{ \mathlarger { \frac{I_{2f}\big|_{p^2}^{2n} + I_{3f}\big|_{p^3}^{2n}}{2n} } } &\text{for} \ 3 \leq log_p(2n) < 4 \\
      \mathlarger{ \mathlarger { \frac{I_{2f}\big|_{p^2}^{2n}} {2n} } } &\text{for} \ 2 \leq log_p(2n)<3  \\
      0 &\text{for}\ log_p(2n) < 2\\ 
    \end{cases}
    \label{eq: alpha}
\end{equation}

\begin{figure}[H]
    \centering
    \includegraphics[width=\linewidth]{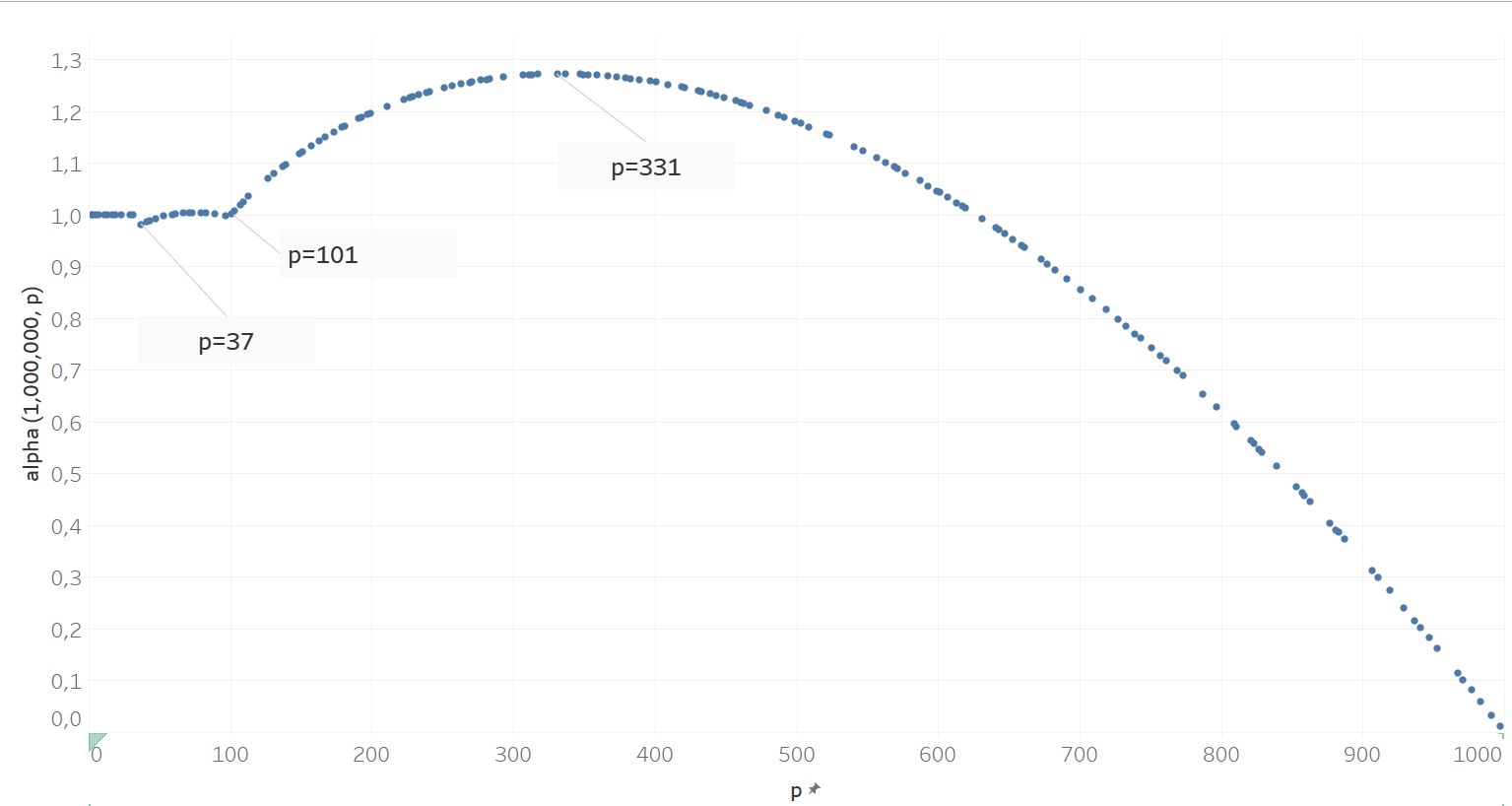}
    \caption{Alpha for all primes in the Pen for 2n = 1 Million}
    \label{fig:rho 1M}
\end{figure}

Fig. \ref{fig:rho 1M} shows $\alpha(1,000,000,p)$. For small primes the value of $\alpha$ is one as $log_p(2n) > 4$. The first prime where $log_p(2n) < 4$ is 37 and the value of $\alpha$ is just below one, as you also expect from Fig. \ref{fig:RPF actual 31}. From $p=101$  onward only the two factor integral is evaluated as $log_{101}(1,000,000) < 3$. As $p$ further increase, the value of $\alpha$ increases as most of the two factor $TRPF$ is larger than one, Fig. \ref{fig:RPF 2factor multiple}. $\alpha$ continues to increases up to $p=331$, then the fact that $TRPF$ is zero between one and $p^2$ causes a decline in the value of $\alpha$.  

Combining Eq. (\ref{eq:EGP}) and Eq. (\ref{eq: alpha}) leads to an improved estimate of the number of GP ($IGP(2n)$). 

\begin{equation}
\begin{split}
\ IGP(2n) & =\frac{n}{2} \ \prod_{p=3}^{p<\sqrt{2n}} \frac{f(2n,p)}{p} \\
    & \text{   where   }  f(2n,p) \begin{cases}
      p-\alpha(n,p),&\text{if}\ mod(2n,p)=0 \\
      p-2\alpha(2n,p),&\text{if}\ mod(2n,p) \ne 0\\
    \end{cases} \\
    \label{eq: IGP}
\end{split}    
\end{equation}

The improvement is due to the fact that the formula now takes into account the variations in $TRPF$ for each prime in the Pen. The fraction that is eliminated is $2\times \alpha(2n,p) /p$, if $p$ is not a divisor of $2n$. If $p$ is a divisor of $2n$ then you need to evaluate $\alpha$ up to $n$, as on the inbound stretch of the racetrack there are no additional multiples of $p$ that eliminate a pair. 

Figure \ref{EGP, IGP and GP} shows the initial estimate (EGP) (blue), the improved estimate (IGP) (red) and actual number of GP (grey) for the band $B_2$. The ICP is within the band of actual number of GP.
The IGP is lower than the EGP due to fact that $\alpha$ is higher than one for most of the primes in the Pen and it is also higher when there are more remaining pairs. The impact of $\alpha$ is higher for smaller primes.

\begin{figure}[H]
    \centering
    \includegraphics[width=\linewidth]{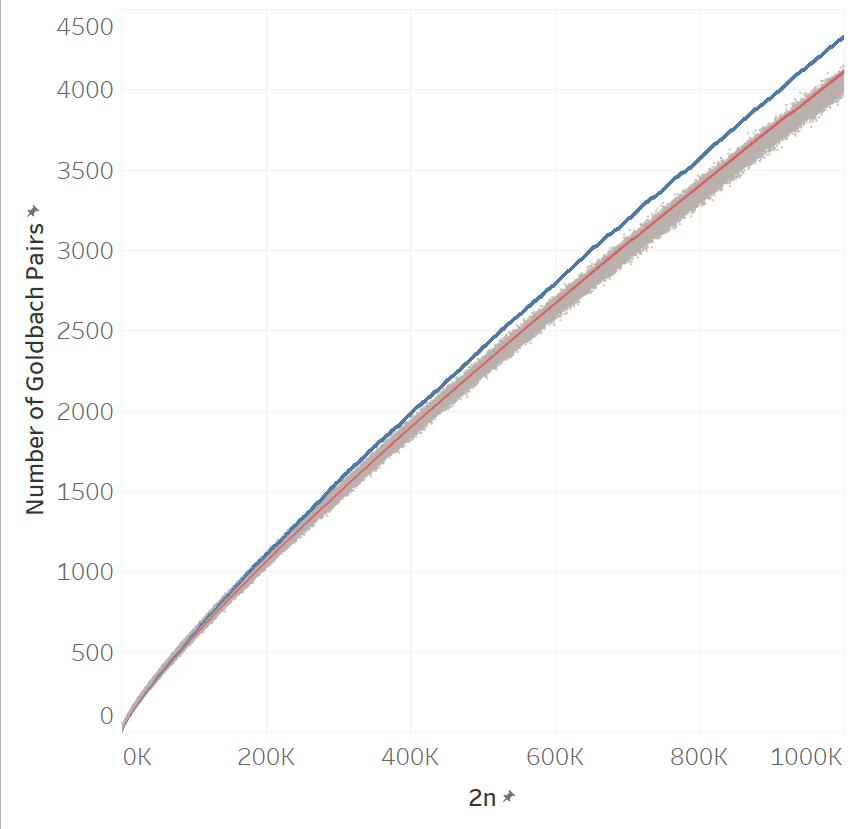}
    \caption{EGP (blue), IGP (red) and actual GP (grey) for Band $B_2$ up to 2n from 2 ... 1 million}
    \label{EGP, IGP and GP}
\end{figure}

\newpage
\section{Road to a proof of the Goldbach Conjecture}

The above is not a formal proof of the Goldbach Conjecture but it does provide a road map on how such a proof might be structured. The key steps are:

\begin{itemize}
  \item Use the racetrack analogy to group $2n$ in Bands based on the divisor primes in the Pen
  \item Band $B_2$ has the lowest estimate for the number of GP and also the lowest actual number of GP. Once the Conjecture is proven for $B_2$ then extending the proof to all 2n should be relatively simple
  \item Show that the function that estimates the number of GP for Band $B_2$ is increasing as $2n$ increases
  \item Further improve the estimate using amongst other the $TRPF$ to come to an unbiased estimator of the number of GP
  \item Quantify the remaining error term between the actual GP and the IGP. Currently the IGP on average still overestimates the actual GP but the bandwidth of the error term for band $B_2$ is small, it seems to be around $2Li(\sqrt{2n})$ 
\end{itemize}

Based on the above steps you would end up with a function that estimates the GP, that is ever increasing and has a error band which is much smaller then the function value. That would proof the Goldbach conjecture that for every even number there is at least one GP.
    
For the third step, the starting point is Eq. {\ref{eq:EGP}} to estimate the number of GP for Band  $B_2$, which simplifies to 

\begin{equation}
 \ EGP_{B2}(2n)=\frac{n}{2} \ \prod_{p=3}^{p<\sqrt{2n}} \frac{p-2}{p} 
 \label{eq: EGP B2}
\end{equation}

\noindent
Rearranging and applying Merten's third theorem leads to 

\begin{equation}
\begin{split}
\ EGP_{B2}(2n)&=\frac{n}{ln^2(\sqrt{2n})}\ ln(\sqrt{2n})\prod_{p=3}^{p<\sqrt{2n}} \frac{p-2}{p-1} \ ln(\sqrt{2n})\prod_{p=2}^{p<\sqrt{2n}} \frac{p-1}{p} \\
\ &=2\frac{2n}{ln^2(2n)}\ c \ e^{-\gamma} \  where\  c \approx 0.74123... \\
\ &= C\frac{2n}{ln^2(2n)} where\  C \approx 0.83235...
\end{split}
\label{eq: EGP B2 rearranged}
\end{equation}

From equation (\ref{eq: EGP B2 rearranged}) it follows that the EGP for $B_2$ is continuously increasing as 2n increases.

For the fourth step the estimation of number of GP needs to be improved. The EGP slight overestimates the number of GP and can be improved by using the $TRPF$. The IGP is close to the actual number of GP but is still not an unbiased estimator of the number of GP. It still seems to slightly overestimate the number of GP for values of $2n$ close to 1 million. There are a number of steps in the calculation of the estimate which could lead to the bias. Amongst others these are:

\begin{itemize}
  \item the assumption in the calculation of $\alpha$ in Eq. ($\ref{eq: alpha}$) that for $4 \leq log_p(2n)$  $\alpha = 1$
  \item the estimate uses $Li(x)$, and $Li(x)$ is slightly larger then the actual number of primes smaller than $x$
  \item the estimate is based on limits as $p \to \infty$ (Eq. (\ref{eq:RPF for 2 factor multiples }) and (\ref{eq: integral 3f})). This may not be correct for smaller $p$ as is also shown by the difference in Fig. \ref{fig:RPF actual 31} and \ref{fig:RPF ln distribution}
  \item the actual distribution of primes is random, so an exact function is not possible.
\end{itemize}

For the final step in the proof, the error terms needs to be quantified and bounded. The actual distribution of primes is random which inevitably leads to an error term in the estimate. Fig. \ref{EGP, IGP and GP} shows that the bandwidth of the actual GP for Band $B_2$ is quite small. One of the key benefits of the road map is that it greatly reduces the error term if you look at individual bands instead of all $2n$. The overall bandwidth in Fig. \ref{fig:bands} is much bigger than the bandwidth of $B_2$ in Fig. \ref{EGP, IGP and GP}.

\medskip
A rigorous proof is required for all the steps on the road map but it may lead to a proof of the Goldbach conjecture.

\newpage
\bibliography{sample}

\end{document}